\documentclass[a4paper]{article}

\usepackage{
amsmath,
amsthm,
amscd,
amssymb,
}
\usepackage{listings}
\usepackage[
            left=3cm,right=3cm,top=3cm,bottom=3cm,%
            ]{geometry}
\usepackage{wasysym}
\usepackage{comment}
\usepackage{tikz}
\usetikzlibrary{positioning,arrows,calc,decorations.pathmorphing}
\usepackage{xspace}

\usepackage{graphicx}

\usepackage{thm-restate}

\usepackage{url}

\theoremstyle{plain}

\newtheorem{theorem}{Theorem}

\newtheorem{question}{Question}

\newtheoremstyle{derp}
{3pt}
{3pt}
{}
{}
{\upshape}
{:}
{.5em}
{}
\theoremstyle{derp}

\newcommand{\R}{\mathbb{R}}

\newcommand{\Z}{\mathbb{Z}}

\newcommand{\N}{\mathbb{N}}

\newcommand{\MCG}{\mathrm{MCG}}

\title{Thompson's $V$ in MCG of mixing SFT by PW-linear homeos}

\author{
Ville Salo \\
vosalo@utu.fi
}

\begin{document}
\maketitle

\begin{abstract}
In a recent paper, we showed that groups admitting ``veelike actions'' on a finite language embed in mapping class groups of certain two-sided subshifts. In this note, we illustrate this theorem for the embedding of Thompson's $V$ by exhibiting the piecewise linear local rules for the embedding. These turn out to split the embedding into the mapping class group, showing that $V$ even embeds in the mapping class group by homeomorphisms.
\end{abstract}

\section{Introduction}

The mapping class group of a subshift was introduced in \cite{Bo02a}, and studied extensively in \cite{BoCh17} (also see \cite{ScYa21}). By embedding suitable groups as subgroups, it was in particular shown that this group is non-amenable and non-residually finite. It was asked in \cite[Question~6.3]{BoCh17} whether this group is sofic. We showed in \cite{Sa21} that solving this question would necessarily solve an open problem, by embedding Thompson's $V$ (and also the Brin-Thompson $2V$). In particular the following theorem was obtained.

\begin{theorem}
\label{thm:V}
Thompson's group $V$ embeds in the mapping class group of the vertex shift defined by the matrix $\left[\begin{smallmatrix}
1 & 1 & 0 \\
1 & 1 & 1 \\
1 & 1 & 1 \\
\end{smallmatrix}\right]$.
\end{theorem}

This is a technical note, where we give a more elaborate discussion of this example: we give explicit piecewise linear local rules implementing the generating set for $V$ given in \cite{BlQu17}, show some example computations and check the defining relations. It turns out that the total flow distortion of each defining relation is trivial, and we obtain the following stronger theorem.

\begin{theorem}
\label{thm:strongV}
The embedding constructed in the proof of Theorem~\ref{thm:V} can be realized with an action of piecewise linear homeomorphisms.
\end{theorem} 

We mean that the map from the group of orientation-preserving homeomorphisms on the mapping torus, defined by quotienting by isotopy, splits on the subgroup we define, by a homomorphism with only piecewise linear homeomorphisms in its image. We are not aware of a theoretical justification for this fact, we found it experimentally and verified it by computation. The theorem is somewhat reminiscent of the Nielsen embedding problem \cite{Ni32,Ke83}.

We assume familiarity with symbolic dynamics \cite{LiMa95}, mapping class groups of SFTs \cite{BoCh17}, and Thompson's group $V$ \cite{BlQu17}. As the present note is a technical supplement to \cite{Sa21}, we use notation from \cite{Sa21} and recommend that the reader starts there. The minor difference is that, to be in line with \cite{BlQu17}, groups act from the right in this note, while they act on the left in \cite{Sa21}.

\section{A concrete implementation of $V$}

Recall that we consider the mapping class group of the vertex shift $X$ defined by the matrix $\left[\begin{smallmatrix}
1 & 1 & 0 \\
1 & 1 & 1 \\
1 & 1 & 1 \\
\end{smallmatrix}\right]$ where the dimensions are indexed by the symbols $0,1,2$ in this order.

We very briefly recall the basic idea of the embedding in \cite{Sa21}. The leftmost point on an interval representing a $2$ in a point of the mapping class torus is called an \emph{anchor}. We interpret the continuation to the right from an anchor as a point in Cantor space $\{0,1\}^\N$: If the mapping torus' natural flow never reaches another anchor, then we directly interpret the continuation as a point of Cantor space, and if we run into an anchor, then (because $02$ is forbidden) the continuation is in a natural correspondence with some finite-support point in Cantor space, and we act on this finite-support configuration. By stretching the flow suitably, we can make sure that far to the right of an anchor the content of the configuration is fixed, and we obtain an embedding of $V$ into the mapping class group.

The mapping class group elements that we will use are then of the following form: In some subwords $uvw$ appearing in the configuration, we replace $v$ by another word $v'$, and stretch the flow linearly so that the time it takes to flow over $v$ is the same as flowing over $v'$ in the image. If in every configuration $x \in X$, every cell is the $v$-piece of exactly one such word $uvw$, then such a rule describes an element of the mapping class group (see \cite{BoCh17} and its references for another approach to local rules). To describe explicit elements of $\MCG(X)$, it then suffices to list the mappings $(u,v,w) \mapsto v'$, and of course we only list words that actually appear in the vertex shift $X$. We write such a mapping briefly as $u(v)w:v'$, and the set of such mappings is called the \emph{local rule}.

Let $a, b, c$ be the generators of $V$ from \cite{BlQu17}. These elements respectively perform the \emph{prefix-permutations} $(00 \;\; 01), (01 \;\; 10 \;\; 11), (00 \;\; 1)$ on the prefix when applied to $x \in \{0,1\}^\N$. For $a$ we pick the following local rule, where $\alpha,\beta,\gamma$ range over $\{0,1\}$:
\[ \alpha\beta(\gamma):\gamma, (2)2:201, (200):201, (201)\alpha:200, (201)2:2, (21)2:21, (21\alpha):21\alpha. \]
For $b$ we pick the following local rule, where $\alpha,\beta,\gamma$ range over $\{0,1\}$:
\[ \alpha\beta(\gamma):\gamma, (2)2:2, (200):200, (201)\alpha:210, \]
\[ (201)2:21, (21)2:211, (210)\alpha:211, (211):201. \]
For $c$ we pick the following local rule, where $\alpha,\beta$ range over $\{0,1\}$:
\[ \alpha0(\beta):\beta, 1(\alpha):\alpha, (2)2:21, (200):21, (201):201, (21)\alpha:200, (21)2:2. \]

The idea behind these choices is that we ``fix the anchors'', and to the right of them, we look for the prefix we know how to rewrite. If this would rewrite the entire word visible (or we do not see enough bits to apply a rewrite), then we use a different logic, and we interpret the continuation instead as a finite-support configuration followed by infinitely many zeros, and in this case we apply the prefix-permutation to this finite-support configuration, remove the trailing zeroes and write the resulting finite word in place of the rewritten word. The rewrite is always performed with a constant slope. This is perhaps best internalized by looking at the spacetime diagrams in the following section and in Appendix~\ref{sec:LookHomeo}. An abstract version of this logic (for a general prefix-rewriting bijection) is implemented as part of the program in Appendix~\ref{sec:CheckHomeo} (this is the function {\color{blue}\texttt{apply}}).

Now, to see that this is a representation by homeomorphisms, it suffices to check that every trivial element is not just isotopic to the identity, but is actually the identity map. It suffices to check the relations of any finite presentation. We check the relations
\[ aa, bbb, cc, abababab, cacaca, cabbabacabbabacbcababbacababba, acbcbabbcbbcbcabcbbcabbacbbcbcbabb, \]
\[ abbcbcabbabbcbcbbabbcbbcbabcbbcabb, cabbcbbcbacabacbcbbcabbcabcbbcbbacbacbcbbcabb. \]
These relations are essentially (2.4) in \cite{BlQu17}; the only difference is that we have removed inverses using $a^{-1} = a$, $b^{-1} = b^2$ and $c^{-1} = c$, and have added the third relation $cc$ so that the last substitution is safe.

Observe that indeed all the generators fix the anchors, so the same is true for the compositions. Thus we only need to analyze the action of each relation $g$ on the list above on words of the form $2w$, where $w \in \{0,1\}^*$ does not end in $0$. More precisely, we can look at their actions on the points of the mapping class torus corresponding to the bi-infinite words $(2w)^\Z$. This reduces the problem to checking the distortion of the action of relations on a countable set of words. To reduce to a finite set of words, observe that each relation will only look a finite distance into a long word $w \in \{0,1\}^*$. Thus, it is enough to check the distortion of the relations when applied to $2w$ for longer and longer words, until the last symbol is not actually touched.

Code for checking this is included in Appendix~\ref{sec:CheckHomeo}. Diagrams for checking lack of distortion by visual inspection are in Appendix~\ref{sec:LookHomeo}. After inspecting and running the code, or carefully scrutinizing these diagrams, we conclude that Theorem~\ref{thm:strongV} holds.

\section{Spacetime diagrams}
\label{sec:Spacetime}

In this section, we show some pictures of what the action of $V$ looks like, with the local rules chosen in the previous section. This section serves no precise mathematical purpose, but the pictures may help understand the construction, and we found Theorem~\ref{thm:strongV} by looking at them and observing that the elements that the theory guarantees are isotopic to the identity actually act trivially on uniformly random finite words with high probability (we first assumed this was a bug in our program). 

We use the convention that a configuration of the mapping class group is shown on the top row, and below it we show the successive images when a sequence of generators $a,b,c$ are applied. We distort the flow in the image configurations for each partial application $g$, so that $t + (x \cdot g) = (t + x) \cdot g$ holds (where $t+y$ denotes the $\R$-flow by $t \in \R$); we refer to this as \emph{cocycle distortion}. We use the following bitmaps to represent the numbers $0, 1, 2$:
\begin{center}
\includegraphics[scale=0.25]{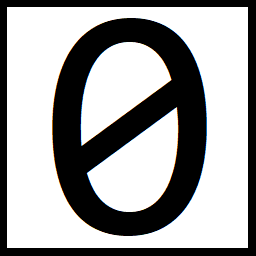} \;\;\;\;\;\; \includegraphics[scale=0.25]{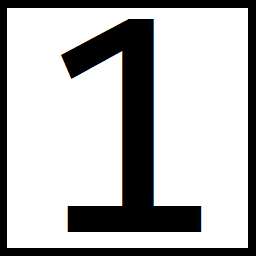} \;\;\;\;\;\;\; \includegraphics[scale=0.25]{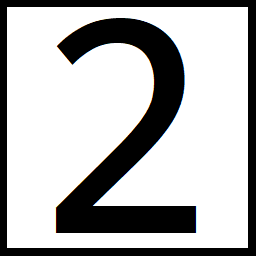}
\end{center}
\noindent Cocycle distortion shows up as literal distortion of these bitmap images. Note that the actual images can be read off by undoing the cocycle distortion (outwards from the chosen origin), by undistorting the bitmaps.

These can be thought of as \emph{spacetime diagrams} for a mapping class group element (more precisely, for the local rule of one), analogously to spacetime diagrams as defined in cellular automaton theory -- indeed reversible cellular automata give elements of mapping class groups, and with the natural choice of local rule, their classical spacetime diagram is the same as that of the corresponding mapping class group element. Note that the natural flow on the mapping torus should be thought of as the flow of \emph{space}, and the mapping class group element gives the (discrete) flow of time.

The following shows the spacetime diagram for the action of $cbcabb$ (composing left to right, thus top to bottom).
\begin{center} 
\begin{tikzpicture}
\node [anchor=south west] (label) at (0,0) {
\includegraphics[scale=0.3,trim={3.85cm 6.5cm 6.45cm 2.5cm},clip]{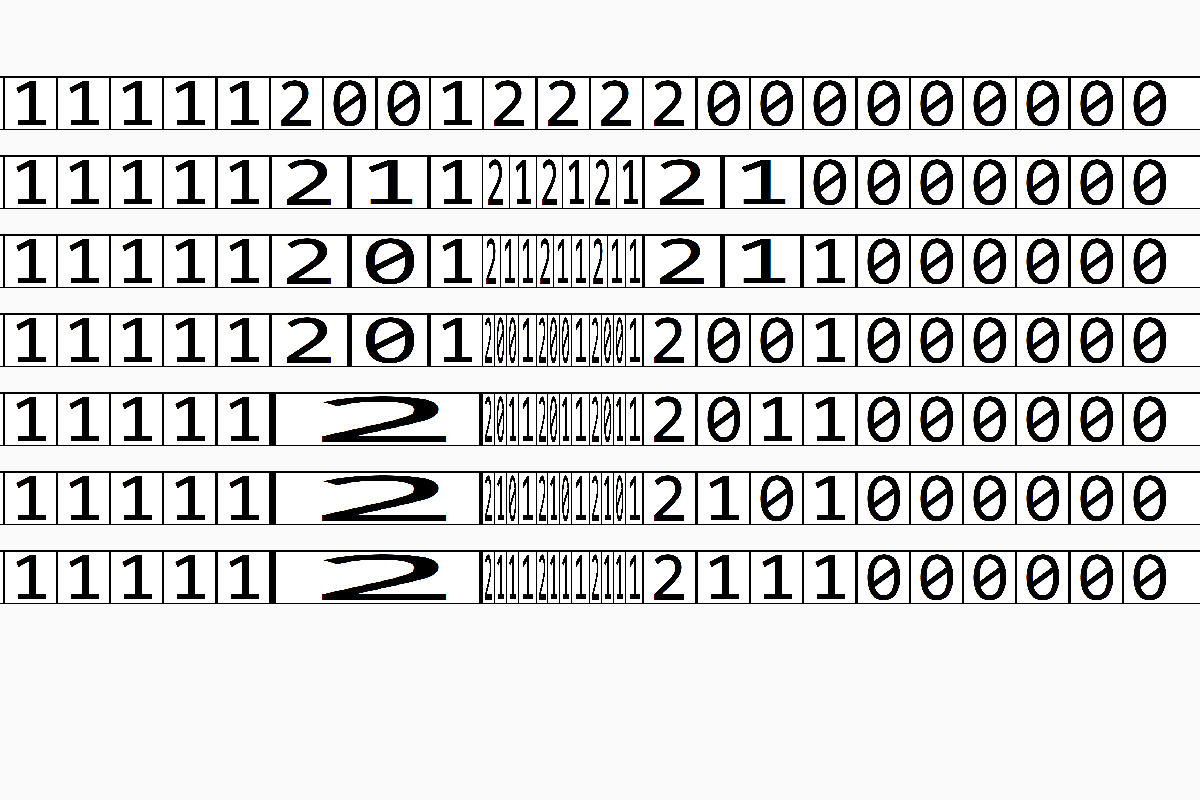}};
\draw (0,4.8) edge[bend left,<-] node[left] {$c$} (0,5.4);
\draw (0,3.97) edge[bend left,<-] node[left] {$b$} (0,4.57);
\draw (0,3.14) edge[bend left,<-] node[left] {$c$} (0,3.74);
\draw (0,2.31) edge[bend left,<-] node[left] {$a$} (0,2.91);
\draw (0,1.48) edge[bend left,<-] node[left] {$b$} (0,2.08);
\draw (0,0.65) edge[bend left,<-] node[left] {$b$} (0,1.25);
\end{tikzpicture}
\end{center}

In the suffix $2000000...$, we simulate the natural action of $V$ with only minor distortion in the flow. The piece $2001$ in the pattern $20012$ represents the infinite word $x = 2001000...$, and $x \cdot cbcabb= 2000...$, so this pattern disappears in the application of $cbcabb$ (thus we shrink  this part by a factor of four, i.e.\ the cocycle moves quickly in this part). The first symbol $2$ in a pattern $22$ represents the configuration $y = 2000...$ and $y \cdot cbcabb = 2111000...$, so the run of $2$s becomes a run of $2111$s (thus we stretch them by a factor of four, i.e.\ the cocycle moves slowly in this part).

By Theorem~\ref{thm:strongV}, any composition of $a,b,c$ which evaluates to identity in $V$, also gives the identity in the mapping class group. It is, however, possible for an element to fix the flow orbit of a configuration while acting non-trivially on it, even on the configurations $...00012000...$ where we simulate the natural action of $V$. An example of this is the following computation corresponding to $0^\N \cdot cbba = 0^\N$.

\begin{center}
\includegraphics[scale=0.18,trim={0 1.3cm 0 1.2cm},clip]{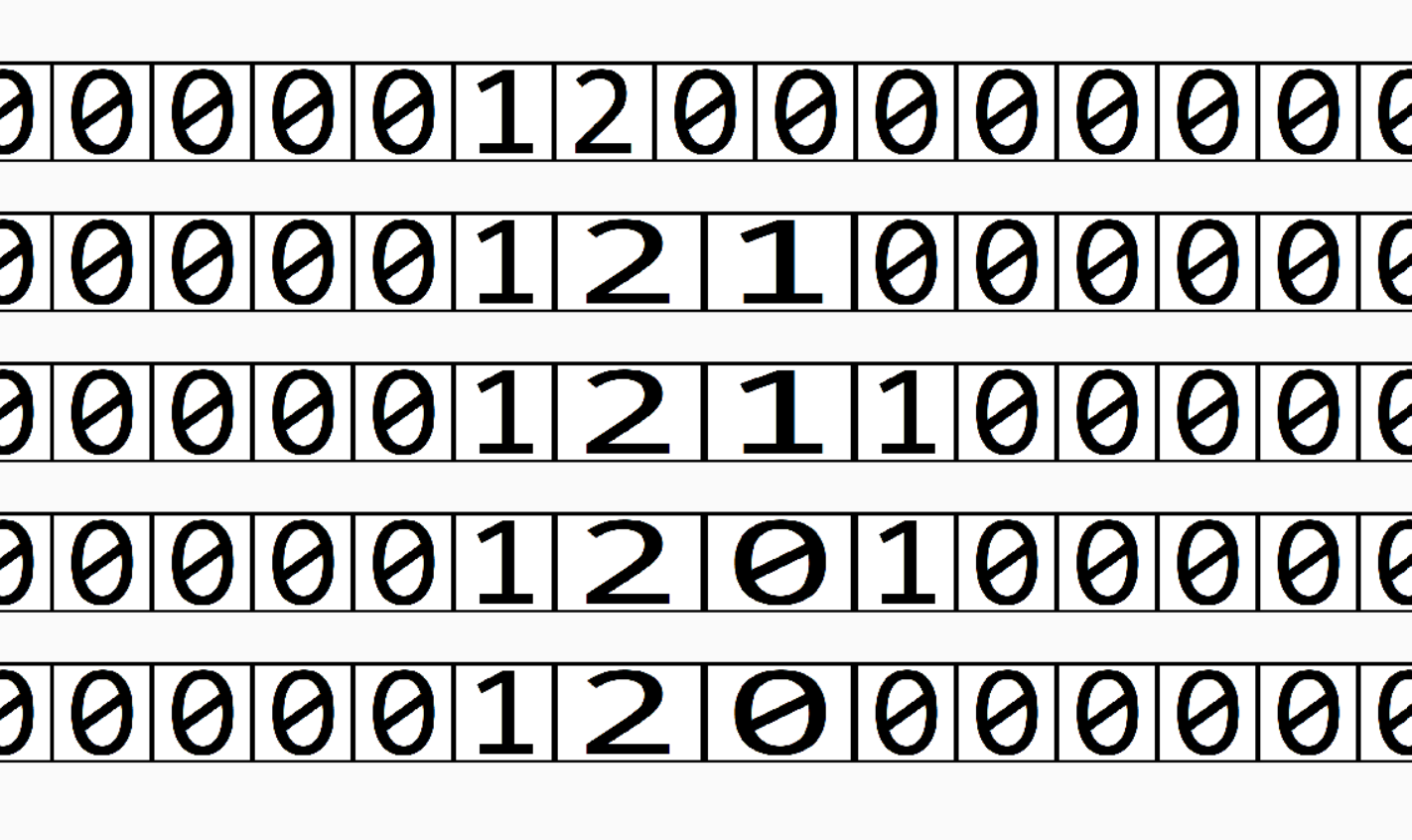}
\end{center}

The explanation is that the natural action of the (non-trivial) $V$-element $cbba$ does not fix the ``germ'' of $0^\N$; intuitively, the cocycle is actually shifting all the zeroes to the left. In the mapping class group, the zeroes far from the ``origin'' $2$ cannot possibly know they are being shifted, meaning we must implement the movement by dilating the prefix of the configuration. This is illustrated in the following figure, which shows the same action with more context.

\begin{center}
\includegraphics[scale=0.3,trim={0 1.4cm 0 1.2cm},clip]{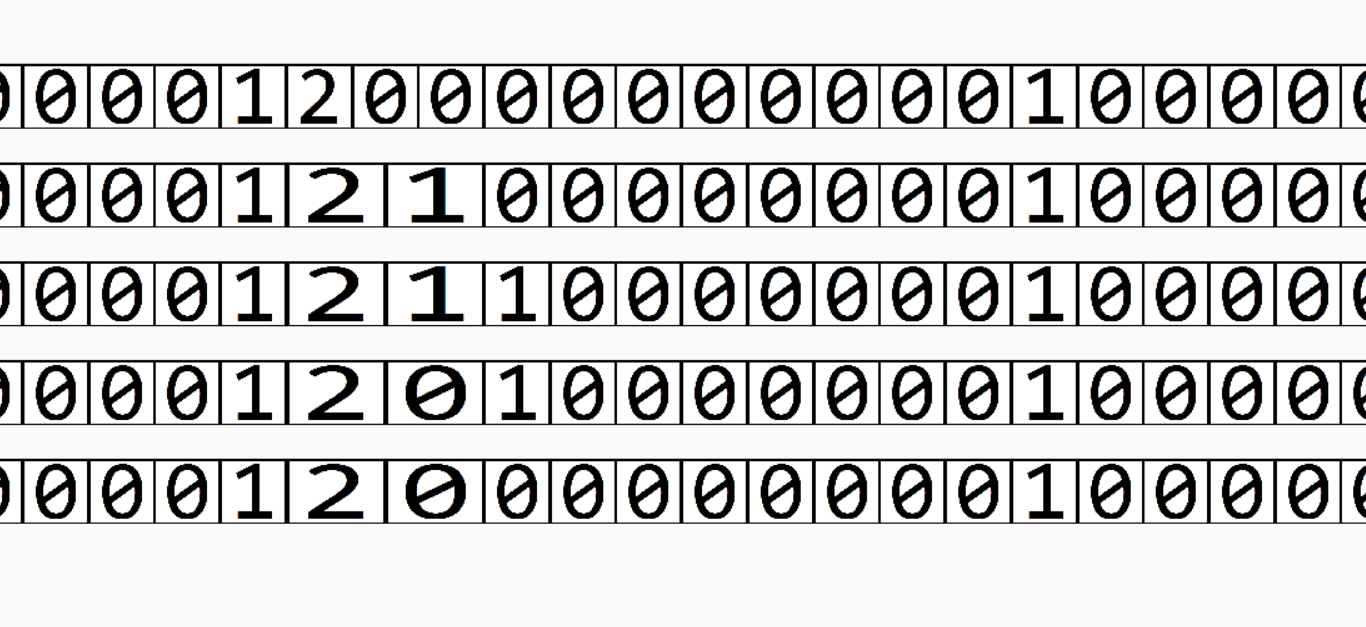}
\end{center}

A piece of a spacetime diagram for the reverse of the penultimate relation, namely
\[ bbacbbcbabcbbcbbabbcbcbbabbacbcbba \]
(which is also a relation) is shown on a random configuration in the following figure.
\begin{center}
\includegraphics[scale=0.75,trim={4cm 0 10cm 0},clip]{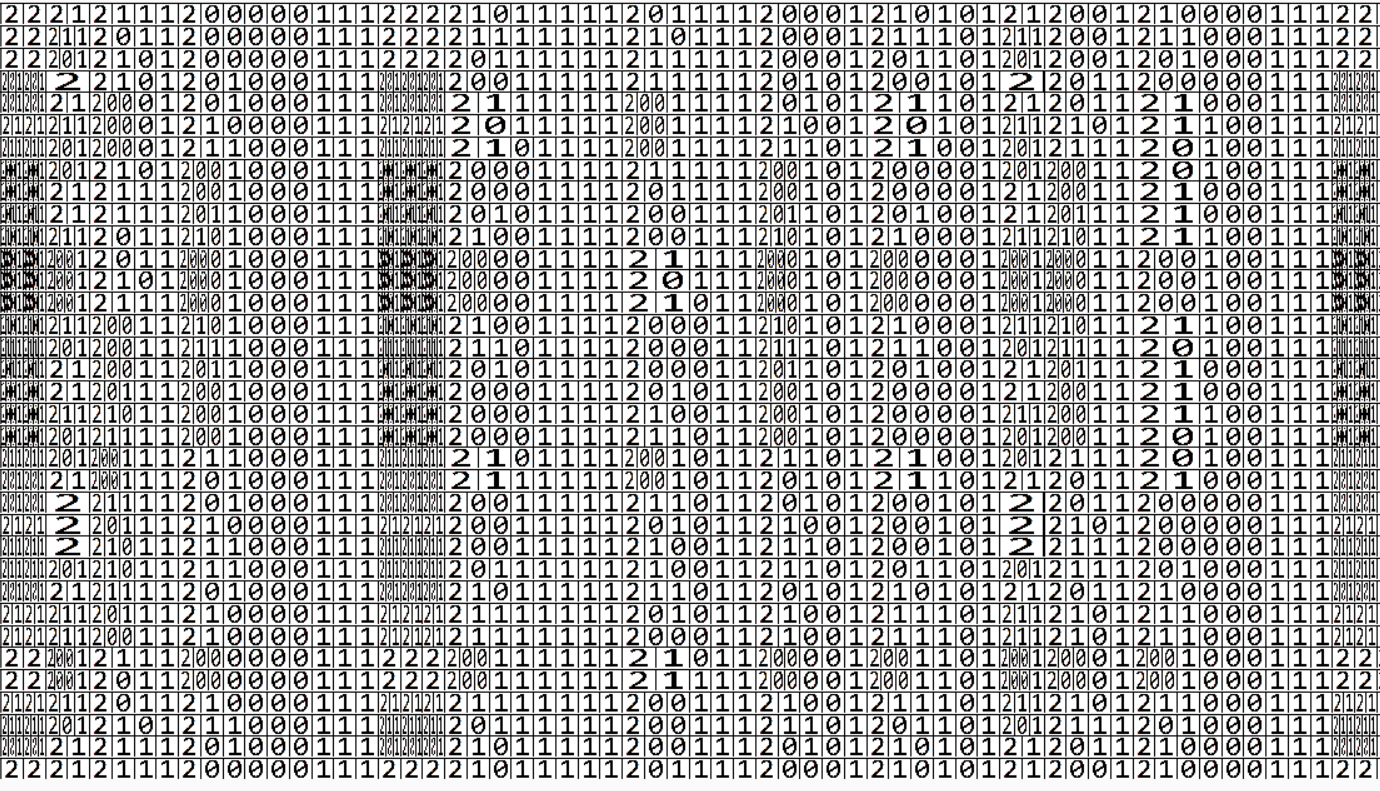}
\end{center}

\section{Discussion}

We note that even having the generators $a,b,c$ act the same way on Cantor space but changing their distortion can break the embedding, for example preserving $a$ and $b$ and replacing $c$ by $000 \mapsto 10, 001 \mapsto 11, 1 \mapsto 00, 01 \mapsto 01$ breaks the embedding of $V$: the action of $cc$ distorts the configuration $(211)^\Z$, thus even breaks the embedding of $\Z_2 \cong \langle c \rangle$. One possible explanation is that our rule for turning a prefix permutation into an MCG element is not the most natural one.

It is open whether the entire mapping class group splits by homeomorphisms in the same sense as our embedding of $V$ does, see \cite[Question~2.1]{BoCh17}. It is also open whether the (Bowen-Franks kernel of the) mapping class group itself is finitely-generated \cite[Question~3.10]{BoCh17} or finitely-presented. If it turns out to be, one could imagine using a similar computational approach to splitting it, although admittedly it seems likely that proving it is finitely presented (if it is) is more difficult than finding a split (if there is one).

It seems unlikely that the present computational approach is really needed to prove Theorem~\ref{thm:strongV}, and it seems likely that there exists a simple explanation. In particular, a more careful analysis of the deduction of this presentation in \cite{BlQu17} might explain this phenomenon.

Because of the point raised in the previous paragraph, and the fact carefully checking the embedding of $V$ was already lot of work, we have not carefully looked at our representation of $2V$ in the mapping class group of a mixing SFT (defined in \cite{Sa21}), to see if it also splits. In theory, the same approach could work, since this group is also finitely presented.

\begin{question}
Does the embedding of $2V$ defined in \cite{Sa21} split by homeomorphisms? If not, does $2V$ admit a split embedding into the mapping class group of a mixing SFT?
\end{question}

\section*{Acknowledgements}

We thank Matt Brin for pointing out the similarity to the Nielsen embedding problem and Laurent Bartholdi for his comments. The author was supported by Academy of Finland project 2608073211.

\bibliographystyle{plain}
\bibliography{../../../bib/bib}{}

\newpage

\appendix

\section{Diagrams for checking lack of distortion}
\label{sec:LookHomeo}

In this section, we include diagrams for checking that there is no distortion in the action of the generators, by visual inspection. We note that these diagrams were produced for completeness, and because they technically constitute a human readable proof of lack of distortion. We have \emph{not} verified all of these diagrams by hand, although we have checked a large sample. 

Rather than reading the exhaustive listing, checking the Python program in Appendix~\ref{sec:CheckHomeo} (and running it) is our recommended way of checking this proof. Alternatively, checking some of these diagrams or playing with this action by hand, one is quickly lead to the warm feeling that everything is going to be all right, and the lack of distortion is actually obvious because of the perfect choice of generators. We have not succeeded in translating this feeling into a proof.

A few words about these diagrams. We use essentially the same conventions as in Section~\ref{sec:Spacetime}, but we flip the axes (the action then goes left-to-right) so that at least one verification fits on a single row. We also use diagrams rather than symbols, and the vertical area between two ticks (short horizontal lines) represents a symbol. A box marks a $1$, a line marks a $0$, a triple line marks the initial $2$, and a bowtie marks an arbitrary continuation of a string of bits (which we are not allowed to touch).

A few (57) discontinuities appear inside tiles, and we mark discontinuities with a gray tick. In the presence of discontinuities, the diagrams do not give full information about the slopes, but at almost all discontinuities it is clear from looking at the diagram that the slope cancels out after a few steps, since the same type of stretching that introduced it is almost immediately used to remove it. The only exceptions are the $bcba$ and $abcb$ fragments of the longest (last) relation applied to the input word $2011$. Here, if one mentally stretches the lowest segment by $3/2$, the cancellations become clear.

We go through words not ending in $0$ in increasing length and in lexicographic order within each length, and then repeat the check at the last length considered (now including words ending in $0$), with the added bowtie.

\input{relation_checks}

\newpage
\section{Code for checking lack of distortion}
\label{sec:CheckHomeo}

The following Python script checks that the relations indeed perform the identity homeomorphism. We mention some implementation details: checking words $w$ up to the possible radius of the relations is tricky, since the longest relation has length $45$. However, the relations do not actually look deeper than four symbols into the word (including the initial $2$ symbol); for the first four relations, three symbols suffice. We can check them up to this length, and also check that we never look deeper, by adding a (meaningless) $3$-symbol at the end of the longest words, and making sure that if this symbol is seen by the local rules, an error is thrown. A (meaningful) $2$-symbol is used in the beginning, as we are also distorting the length of this initial symbol in each segment beginning with $2$. For additional robustness we check words up to length $11$.

The way we check for distortion computationally is that we reprent a word as a list of tuples of the form $(s,a,b,c,d)$ where $s$ is a symbol, $[a,b]$ is a subinterval in the length-$1$ $s$-tile that the piece spans, and $[c, d]$ is a physical interval where this piece is stretched. Every time the $[a,b]$ interval covers $[0,1]$, we change the symbol. This represents what the action of the element would do to a word between two hypothetical $2$s.

\begin{lstlisting}[language=Python]
from fractions import Fraction
div = Fraction # the Fraction constructor also performs exact division
from functools import reduce

a = ["00", "01"]
b = ["01", "10", "11"]
c = ["1", "00"]

NOT_RELATION = "not a relation"
DISTORTED = "distorted"
ERROR = "error"
OK = "ok"

def concat(a,b):return a+b

def binary_words(n):
    if n == 0:
        yield ""
    else:
        for w in binary_words(n-1):
            yield w + "0"
            yield w + "1"

def words_not_ending_in_0(n):
    if n == 0:
        yield ""
    else:
        for w in binary_words(n-1):
            yield w + "1"

"""
Turn a word to a list of tiles.
A tuple (s, a, b, c, d) means the "abstract" [a, b] piece of an s-tile
is linearly stretched on the "physical" [c, d] interval.
"""
def word_to_tiles(word):
    tiles = []
    for i,s in enumerate(word):
        tiles.append((s,0,1,i,i+1))
    return tiles

"""
Remove meaningless cuts in abstract tiles
i.e. cuts where the derivative is continuous.
"""
def simplify(pieces):
    pieces = pieces[:] # always copy for safety
    changed = True
    i = 0
    while i < len(pieces)-1:
        # if not inside an abstract tile, continue
        if pieces[i][2] == 1:
            i += 1
            continue
        # check derivatives
        deriv1 = div(pieces[i][4]-pieces[i][3], pieces[i][2]-pieces[i][1])
        deriv2 = div(pieces[i+1][4]-pieces[i+1][3], pieces[i+1][2]-pieces[i+1][1])
        if deriv1 == deriv2:
            pieces[i] = (pieces[i][0], pieces[i][1],
                         pieces[i+1][2], pieces[i][3], pieces[i+1][4])
            del pieces[i+1]
        else:    
            i+=1
    return pieces

"""
Apply a prefix rewrite to a sequence of tiles.
This implements, abstractly, the same logic we used to construct
the MCG embedding of V.
"""
def apply(rewrite, tiles):
    # collect abstract length one pieces
    pieces = []
    piece = []
    for t in tiles:
        piece.append(t)
        if t[2] == 1:
            pieces.append(piece)
            piece = []
    word = []
    for p in pieces:
        word.append(p[0][0])
    # keep the string representation (for no reason)
    word = "".join(word)
    assert word[0] == "2"
    """
    Next, we figure out how many symbols are rewritten and by what.
    The rule is, if we can actually rewrite a prefix, we do that
    and linearly span. If we cannot (including the case where we
    rewrite the entire word), we add zeroes, rewrite a prefix,
    and remove all tailing zeroes. (And the entire word is rewritten.)
    The only variables computed between here and the next
    comment are removed_count and new_prefix (and new_prefix_length).
    """
    word_wo_2 = word[1:]
    result = prefix_rewrite(rewrite, word_wo_2)
    if result == None or result[0] == len(word_wo_2):
        # first rewrite failed, retry with a long zero suffix
        word_0s = word_wo_2 + "0"*max(map(len, rewrite))
        result = prefix_rewrite(rewrite, word_0s)
        if result == None:
            # this means we ran into 3, which should not happen
            return None
        else:
            # we added zeroes, now apply, remove zeroes,
            # and linearly span the entire word
            _, new_prefix, old_suffix = result
            new_prefix = "2" + new_prefix + old_suffix # reinsert 2
            while len(new_prefix) > 0 and new_prefix[-1] == "0":
                new_prefix = new_prefix[:-1]
            removed_count = len(word_wo_2) + 1
            new_prefix_length = len(new_prefix)
        global max_denom        
    else:
        # first rewrite succeeded
        removed_count, new_prefix, old_suffix = result
        removed_count += 1 # reinsert 2
        new_prefix = "2" + new_prefix # reinsert 2
        new_prefix_length = len(new_prefix)
    
    """
    We now construct tiles for new_prefix, rescale them so they
    take up the same physical space as the first removed_count
    pieces. Then we cut them up
    """
    new_pieces_suffix = pieces[removed_count:]
    new_pieces_prefix = slerp_slorp(pieces[:removed_count], new_prefix)
    return simplify(new_pieces_prefix + reduce(concat,new_pieces_suffix,[]))

"""
Compute the barycentric position of a on [b,c] and
give same barycentric point on [d,e].
"""
def barylerp(a,b,c,d,e):
    t = (a - b)/(c - b)
    return d + t*(e - d)

"""
Make sure a piecewise linear specification contains position in its domain.
"""
def insert_keypoint(keypoints, position):
    if position == 0:
        return keypoints
    new_keypoints = {}
    inserted = False
    keys = list(sorted(keypoints))
    for i,k in enumerate(keys):
        if k == position:
            return keypoints
        if not inserted and k > position:
            b = barylerp(position, keys[i-1], keys[i],
                         keypoints[keys[i-1]], keypoints[keys[i]])
            new_keypoints[position] = b
            inserted = True
        new_keypoints[k] = keypoints[k]
    return new_keypoints

"""
Make sure there is a tile cut at position, on the abstract side.
This assumes the "tiles" representation, so (s, a, b, c, d) tuples,
and a, b always run from 0 to 1 indexing positions on abstract tiles,
tile by tile, while c, d run over physical spans.
"""
def insert_cut(tiles, position):
    if int(position) == position:
        return tiles
    new_tiles = []
    tot = 0
    inserted = False
    for t in tiles:
        currleft = tot+t[1]
        currright = tot+t[2]
        # position is already there
        if position == currleft or position == currright: return tiles
        # position to insert is inside this tile
        if not inserted and currleft < position < currright:
            lerppos = barylerp(position-tot,t[1],t[2],t[3],t[3])
            new_tiles.append((t[0],t[1],position-tot,t[3],lerppos))
            new_tiles.append((t[0],position-tot,t[2],lerppos,t[4]))
            inserted = True
        else:
            new_tiles.append(t)
        # abstract tile changes
        if t[2] == 1: tot += 1
    return new_tiles

"""
Change physical spans of tiles to be according to keypoints.
Tiles should be cut at keypoints already.
"""
def change_physics_by_keypoints(tiles,keypoints):
    new_tiles = []
    tot = 0
    for t in tiles:
        new_tiles.append((t[0], t[1], t[2],
                          keypoints[t[1]+tot], keypoints[t[2]+tot]))
        # abstract tile changes
        if t[2] == 1: tot += 1
    return new_tiles

def scale_keypoints_domain(keypoints, r):
    new_keypoints = {}
    for k in keypoints:
        new_keypoints[k*r] = keypoints[k]
    return new_keypoints

"""
Slerp slorping is easily understood in the case when the pieces variable only
contains tiles with length one (e.g. comes form word_to_tiles). In this case,
the word word is turned into tiles which all take up the same physical length,
stretching each letter by the same amount. We simply conjugate this transform
by the piecewise linear distortion of pieces.
"""
def slerp_slorp(pieces, word):
    tiles = word_to_tiles(word)
    keypoints = {}
    tot = 0 # total length of abstract tiles
    for p in pieces:
        for c in p:
            keypoints[c[1] + tot] = c[3]
        tot += c[2]
    assert tot == len(pieces)
    keypoints[tot] = c[4]
    # we are actually changing the piecewise linear
    # curve of word, so scale the keypoints
    keypoints = scale_keypoints_domain(keypoints, div(len(word), tot))
    for i in range(len(word)):
        keypoints = insert_keypoint(keypoints, i)
    for k in keypoints:
        tiles = insert_cut(tiles, k)
    tiles = change_physics_by_keypoints(tiles, keypoints)
    return tiles

def word_on_tiles(tiles):
    word = ""
    for t in tiles:
        if t[2] == 1:
            word += t[0]
    return word

def orbit_is_undistorted(relator, word):
    original_tiles = word_to_tiles(word)
    tiles = word_to_tiles(word)
    for i,gen in enumerate(relator):
        newtiles = apply(gen, tiles)
        if tiles == None:
            return ERROR
        tiles = newtiles
    if original_tiles == tiles:
        return OK
    if word_on_tiles(tiles) == word:
        return DISTORTED
    return NOT_RELATION

"""
Try to rewrite a prefix.
Return how many symbols were removed, what was added and what's left,
or if rewriting fails return None.
"""
def prefix_rewrite(rewrite, word):
    for r in rewrite:
        if word[:len(r)] == r:
            return len(r), rewrite[r], word[len(r):]
    # Failure: we ran into 3 or we ran out of length.
    return None

relations = ["aa", "bbb", "cc", "abababab", "cacaca",
             "cabbabacabbabacbcababbacababba",
             "acbcbabbcbbcbcabcbbcabbacbbcbcbabb",
             "abbcbcabbabbcbcbbabbcbbcbabcbbcabb",
             "cabbcbbcbacabacbcbbcabbcabcbbcbbacbacbcbbcabb"]

gens_to_permus = {"a":{"00":"01", "01":"00", "10":"10", "11":"11"},
           "b":{"00":"00","01":"10", "10":"11", "11":"01"},
           "c":{"00":"1", "1":"00", "01":"01"}}
   
everything_ok = True
not_a_relation = None
distorted_relation = None
depths = [10]*9 # optimal lengths: [2,2,2,2,3,3,3,3,3]

for dep, relation in zip(depths, relations):
    for l in range(dep+2):
        if l <= dep:
            words = words_not_ending_in_0(l)
        else:
            words = binary_words(dep)
        for w in words:
            
            if l == dep+1: w += "3"

            g2b = [gens_to_permus[a] for a in relation]
            result = orbit_is_undistorted(g2b, "2" + w)
            
            if result == NOT_RELATION:
                everything_ok = False
                not_a_relation = relation
                notrelation_witness = w
                break
            if result == DISTORTED:
                everything_ok = False
                distorted_relation = relation
                distortion_witness = w
            if result == ERROR:
                everything_ok = False
                all = lost
            if not everything_ok: break
        if not everything_ok: break
    if not everything_ok: break

if everything_ok:
    print ("Everything is ok.")
else:
    if not_a_relation:
        print ("Very bad:", not_a_relation,
               "is not even a relation on", notrelation_witness)
    else:
        print ("Bad:", distorted_relation,
               "distorts", distortion_witness)


\end{lstlisting}

\end{document}